\documentclass[preprint,11pt]{elsarticle}

\usepackage{amssymb}
\usepackage{amsmath}
\usepackage{amsfonts}
\usepackage{amssymb}
\usepackage{indentfirst,latexsym,bm}
\usepackage{amsthm}
\usepackage{color,xcolor}
\usepackage[all]{xy}
\input{mathrsfs.sty}
\newtheorem{theorem}{Theorem}[section]
\newtheorem{lemma}[theorem]{Lemma}
\newtheorem{corollary}[theorem]{Corollary}

\theoremstyle{definition}
\newtheorem{definition}{Definition}[section]

\theoremstyle{definition}
\newtheorem{example}{Example}[section]

\theoremstyle{remark}
\newtheorem{remark}{Remark}[section]
\theoremstyle{question}

\theoremstyle{problem}

\numberwithin{equation}{section}

\journal{XXX}

\begin{document}

\begin{frontmatter}





\title{Characterizations of the semi-harmonious and harmonious quasi-projection pairs on Hilbert $C^*$-modules}
\author[shnu]{Xiaoyi Tian\corref{cor1}}
\cortext[cor1]{Current address: Guangdong Police College, Guangzhou 510230, PR China}
\ead{tianxytian@163.com}
\author[shnu]{Qingxiang Xu}
\ead{qingxiang\_xu@126.com}
\author[HS]{Chunhong Fu}
\ead{fchlixue@163.com}
\address[shnu]{Department of Mathematics, Shanghai Normal University, Shanghai 200234, PR China}
\address[HS]{Health School Attached to Shanghai University of Medicine $\&$ Health Sciences,   Shanghai 200237, PR China}

\begin{abstract} For each adjointable idempotent $Q$ on a Hilbert $C^*$-module $\mathscr{H}$, a specific projection $m(Q)$ called the matched projection of $Q$ was introduced recently  due to the characterization of the minimum value among all the distances from projections to  $Q$. Inspired by the relationship between $m(Q)$ and $Q$, another term called the quasi-projection pair $(P,Q)$ was also introduced recently, where $P$ is a projection on $\mathscr{H}$ satisfying $Q^*=(2P-I)Q(2P-I)$, in which
$Q^*$ is the adjoint operator of the idempotent $Q$ and $I$ is the identity operator on $\mathscr{H}$. This paper aims to  make
systematical characterizations of the semi-harmonious and harmonious quasi-projection pairs on Hilbert $C^*$-modules, and meanwhile to provide examples illustrating the non-triviality of the associated characterizations.
\end{abstract}

\begin{keyword} Hilbert $C^*$-module; projection; idempotent;  polar decomposition; orthogonal complementarity.
\MSC 46L08; 47A05



\end{keyword}

\end{frontmatter}




\section{Introduction }\label{sec:Intro}

Unless otherwise specified, throughout this paper, $\mathbb{C}$ is the complex field,
For every $C^*$-algebra $\mathfrak{B}$, $M_n(\mathfrak{B})$ denotes the $n\times n$ matrix algebra over $\mathfrak{B}$.
Specifically, $M_n(\mathbb{C})$ stands for the set of $n\times n$ complex matrices.  Let $\mathscr{H}$ be a (right) Hilbert module over a $C^*$-algebra $\mathfrak{A}$ \cite{Lance,MT,Paschke}.
For each adjointable idempotent $Q$ on $\mathscr{H}$, a specific projection $m(Q)$, called the matched projection of $Q$, is introduced recently in \cite{TXF02} due to the characterization of the minimum value among all the distances from projections on $\mathscr{H}$ to  $Q$. Inspired by the relationship between $m(Q)$ and $Q$, another term called the quasi-projection pair $(P,Q)$   is also introduced in \cite{TXF02}
in the sense of $$Q^*=(2P-I)Q(2P-I),$$ in which $P$ is a projection on $\mathscr{H}$, $Q^*$ denotes the adjoint operator of the idempotent $Q$ and $I$ denotes the identity operator on $\mathscr{H}$. In particular,
$\big( m(Q),Q\big)$
is a quasi-projection pair, which is called the matched pair of $Q$.
It is known that an adjointable operator on a Hilbert $C^*$-module
may have no polar decomposition (see e.g.\,\cite[Example~3.15]{LLX}).
In view of this dissimilarity between operators on Hilbert $C^*$-modules and Hilbert spaces,
the semi-harmonious quasi-projection pair $(P,Q)$ is introduced in \cite{TXF04}.
More precisely, a quasi-projection pair $(P,Q)$ is said to be semi-harmonious if both $P(I-Q)$ and $(I-P)Q$ have the polar decompositions \cite[Definition~2.2]{TXF04}.
Some aspects of the semi-harmonious quasi-projection pairs are dealt with in \cite{TXF04}. To get a deeper understanding of the semi-harmonious quasi-projection pair,  it is worthwhile to investigate some equivalent conditions under which a quasi-projection pair turns out to be
semi-harmonious, which is  the main concern of Section~\ref{sec:semi-harmony} in this paper. Some equivalent characterizations for the semi-harmonious quasi-projection pair are derived in Theorem~\ref{thm:simplified condition for semi-qpp}. Specifically,
it is proved therein that the condition for a quasi-projection pair $(P, Q)$ to be  semi-harmonious can be weakened as one of  $P(I-Q)$ and $(I-P)Q$ has the polar decomposition.
Some other descriptions of the semi-harmonious quasi-projection pair are carried out in
Corollary~\ref{cor:1 to 4}, Lemma~\ref{lem:trivial-01} and Theorem~\ref{thm:for counterexample-01}, with special attention paid to the  semi-harmonious matched pair in Theorem~\ref{thm:semi-harmony characterization for matched pair}.  A method  to construct the semi-harmonious quasi-projection pairs is revealed by Theorem~\ref{thm:invariant submodule M}.
In addition, two types of the non-semi harmonious quasi-projection pairs are provided at the end of Section~\ref{sec:semi-harmony}; see Examples~\ref{ex:not semi-harmonious-11} and \ref{prop:counter example-1} for the details.

Another main concern of this paper is to introduce and clarify the harmonious
quasi-projection pair (see Section~\ref{sec:harmony}),  thereby to establish theoretical tools for our further studying the block  matrix representations of the quasi-projection pairs \cite{TXF05}.
A quasi-projection $(P,Q)$ on a Hilbert $C^*$-module $\mathscr{H}$ is said to be harmonious if both
$\overline{\mathcal{R}\big(PQ(I-P)\big)}$ and $\overline{\mathcal{R}\big((I-P)QP\big)}$  are  orthogonally complemented
in $\mathscr{H}$ (see Definition~\ref{defn of harmonious quasi-projection pair}).
Likewise, some characterizations of the general harmonious quasi-projection pair are given in Theorem~\ref{thm:harmony=2 semi-harmony} and  Corollary~\ref{cor:four equivalences}.
The equivalence of a matched pair to be semi-harmonious and harmonious is revealed by  Theorem~\ref{thm:matched pair semi harmony is harmony}.
It is known that every $C^*$-algebra $\mathfrak{B}$ can be regarded as a Hilbert  module over itself in the natural way.
When $\mathfrak{B}$ is unital, some special quasi-projection pairs can be constructed as \eqref{equ:nnew defn of P Q unital Cstar alg} on the induced  matrix algebra
 $M_2(\mathfrak{B})$.   It is interesting to find out a necessary and sufficient condition under which these quasi-projection pairs are harmonious, which is fulfilled in Theorem~\ref{thm:for harmony in unital alg}. Attention has been paid to construct the semi-harmonious quasi-projection pairs that are not harmonious (see Examples~\ref{ex:semi-harmony does not imply harmony-01} and \ref{ex:semi-harmonious not implication}).

The rest of this paper is organized  as follows. In Section \ref{sec:preliminaries}, we provide
some preliminaries which will be used in the sequel, and  some basic knowledge about the quasi-projection pair and the matched projection are given. Sections~\ref{sec:semi-harmony} and \ref{sec:harmony} are devoted to systematical characterizations of
the  semi-harmonious and  harmonious quasi-projection pairs, respectively.

\section{Preliminaries}\label{sec:preliminaries}

Hilbert $C^*$-modules are the natural generalizations of Hilbert spaces by allowing the inner product to take values in a $C^*$-algebra rather than in the complex field. Given Hilbert $C^*$-modules $\mathscr{H}$ and $\mathscr{K}$ over a $C^*$-algebra $\mathfrak{A}$, let $\mathcal{L}(\mathscr{H, K})$ denote the set of all adjointable operators from $\mathscr{H}$ to $\mathscr{K}$.
It is known that every adjointable operator $A\in\mathcal{L}(\mathscr{H, K})$ is a bounded linear operator, which is also $\mathfrak{A}$-linear in the sense that
\begin{equation}\label{equ:keep module operator}A(xa)=A(x)a\quad (x\in \mathscr{H},  a\in\mathfrak{A}).
\end{equation}
Let $\mathcal{R}(A)$, $\mathcal{N}(A)$, $A^*$  and $A|_\mathscr{M}$ denote the range,  the null space and the adjoint operator of $A$, and the restriction of $A$ on a subset $\mathscr{M}$ of $\mathscr{H}$, respectively.
Let  $\overline{\mathcal{R}(A)}$ denote the norm closure of $\mathcal{R}(A)$.
We use the abbreviation $\mathcal{L}(\mathscr{H})$ whenever $\mathscr{H}=\mathscr{K}$. Let $\mathcal{L}(\mathscr{H})_{\mbox{sa}}$ (resp.\,$\mathcal{L}(\mathscr{H})_+$) denote the set of all self-adjoint (resp.\,positive) elements in the $C^*$-algebra $\mathcal{L}(\mathscr{H})$. In the case that $A\in\mathcal{L}(\mathscr{H})_+$, the notation $A\ge 0$ is also used to indicate that $A$ is a positive operator on $\mathscr{H}$.
For each $T\in\mathcal{L}(\mathscr{H})$, let $\sigma(T)$ denote the spectrum of $T$ with respect to $\mathcal{L}(\mathscr{H})$.

In the remainder of this paper, $\mathfrak{A}$ is a $C^*$-algebra, $\mathscr{E}$, $\mathscr{H}$ and $\mathscr{K}$ are Hilbert $\mathfrak{A}$-modules.

Given a unital $C^*$-algebra  $\mathfrak{B}$, its unit is denoted by $I_{\mathfrak{B}}$. When $\mathscr{V}$ is a Hilbert module over a $C^*$-algebra, the unit of $\mathcal{L}(\mathscr{V})$ will be denoted simply by $I_{\mathscr{V}}$ rather than $I_{\mathcal{L}(\mathscr{V})}$. In most cases of this paper,
 the Hilbert $C^*$-module under consideration is represented by $\mathscr{H}$, so the symbol $I_{\mathscr{H}}$ is furthermore simplified to be $I$. In other words, the capital Italic letter $I$ is reserved for the identity operator on the Hilbert $C^*$-module  $\mathscr{H}$.

 An operator $Q\in\mathcal{L}(\mathscr{H})$ is called an idempotent if $Q^2=Q$.
If furthermore $Q$ is self-adjoint, then $Q$ is referred to be a projection.
A closed submodule $\mathscr{M}$ of  $\mathscr{H}$ is referred to be
orthogonally complemented  in $\mathscr{H}$ if $\mathscr{H}=\mathscr{M}+ \mathscr{M}^\bot$, where
$$\mathscr{M}^\bot=\big\{x\in \mathscr{H}:\langle x,y\rangle=0\ \mbox{for every}\ y\in
\mathscr{M}\big\}.$$
In this case, the projection from $\mathscr{H}$ onto $\mathscr{M}$ is denoted by $P_\mathscr{M}$.
Let $v^T$ denote the transpose of a vector $v$ in a partitioned linear space. Put
$$\mathscr{H}\oplus \mathscr{K}=\left\{(h_1, h_2)^T :h_1\in \mathscr{H},  h_2\in \mathscr{K}\right\},$$ which is also a Hilbert $\mathfrak{A}$-module whose $\mathfrak{A}$-valued inner product is given by
$$\left<(x_1, y_1)^T, (x_2, y_2)^T\right>=\big<x_1,x_2\big>+\big<y_1,
y_2\big>$$ for $x_i\in \mathscr{H}, y_i\in \mathscr{K}, i=1,2$.
It is well-known that up to unitary equivalence, every element in $\mathcal{L}(\mathscr{H}\oplus \mathscr{K})$
can be identified with a block matrix operator. For details, see the monograph \cite{Moslehian-} and the references therein.

The polar decomposition constitutes a fundamental tool for handling operators on Hilbert $C^*$-modules. The theory for unbounded regular operators in this setting was developed by Frank and Sharifi \cite{frank}; however, this paper  focuses exclusively on the bounded case. For each $T\in\mathcal{L}(\mathscr{H, K})$, the notation $|T|$ is used to denote the square root of $T^*T$. Hence,  $|T^*|=(TT^*)^\frac12$. From \cite[Lemma~3.9]{LLX} we know that there exists at most a partial isometry
$U\in\mathcal{L}(\mathscr{H, K})$  satisfying
\begin{equation}\label{equ:two conditions of polar decomposition}T=U|T|\quad \mbox{with}\quad U^*U=P_{\overline{\mathcal{R}(T^*)}}.\end{equation}
If such a partial isometry $U$ exists, then the representation  \eqref{equ:two conditions of polar decomposition}
is called the polar decomposition of $T$ \cite[Definition~3.10]{LLX}. The following lemma gives a characterization of the existence of the polar decomposition.

\begin{lemma}\label{lem:polar decomposition of T star}{\rm \cite[Lemma~3.6 and Theorem~3.8]{LLX}} For every $T\in\mathcal{L}(\mathscr{H, K})$, the following statements are equivalent:
\begin{enumerate}
\item[{\rm (i)}] $T$ has the polar decomposition;
\item[{\rm (ii)}]$T^*$ has the polar decomposition;
\item[{\rm (iii)}]$\overline{\mathcal{R}(T)}$ and $\overline{\mathcal{R}(T^*)}$ are orthogonally complemented in $\mathscr{K}$ and $\mathscr{H}$, respectively.
\end{enumerate}
If {\rm (i)--(iii)} are satisfied, and the polar decomposition of $T$ is given by \eqref{equ:two conditions of polar decomposition}, then the polar decomposition of $T^*$ is represented by
\begin{equation*}\label{equ:the polar decomposition of T star-pre stage}T^*=U^*|T^*|\quad \mbox{with}\quad UU^*=P_{\overline{\mathcal{R}(T)}}.
\end{equation*}
\end{lemma}

\begin{remark}
Suppose that $T\in\mathcal{L}(\mathscr{H, K})$ has the polar decomposition represented by \eqref{equ:two conditions of polar decomposition}. For simplicity, henceforth  we just say that $T$ has the polar decomposition $T=U|T|$.
\end{remark}

We  provide some additional lemmas focusing on the norm closures of the operator ranges, which will be used in the sequel.

\begin{lemma}\label{lem:Range closure of TT and T} {\rm\cite[Proposition 3.7]{Lance}}
For every $T\in\mathcal{L}(\mathscr{H, K})$, $\overline {\mathcal{R}(T^*T)}=\overline{ \mathcal{R}(T^*)}$ and $\overline {\mathcal{R}(TT^*)}=\overline{ \mathcal{R}(T)}$.
\end{lemma}

\begin{lemma}\label{lem:Range Closure of T alpha and T}{\rm (\cite[Proposition~2.9]{LLX} and \cite[Lemma 2.2]{Vosough-Moslehian-Xu})} Let $T\in \mathcal{L}(\mathscr{H})_+$. Then for every $\alpha>0$, $\overline{\mathcal{R}(T^{\alpha})}=\overline{\mathcal{R}(T)}$ and $\mathcal{N}(T^{\alpha})=\mathcal{N}(T)$.
\end{lemma}

\begin{lemma}\label{lem:rang characterization-1}{\rm \cite[Proposition~2.7]{LLX}} Let $A\in\mathcal{L}(\mathscr{H, K})$ and $B,C\in\mathcal{L}(\mathscr{E,H})$ be such that $\overline{\mathcal{R}(B)}=\overline{\mathcal{R}(C)}$. Then $\overline{\mathcal{R}(AB)}=\overline{\mathcal{R}(AC)}$.
\end{lemma}

Now, we turn to recall some basic knowledge about the quasi-projection pair  described in \cite{TXF02}.

\begin{definition}\label{defn:quasi-projection pair}\rm{ \cite[Definition~2.2]{TXF02}} An ordered pair $(P,Q)$  is called a quasi-projection pair on $\mathscr{H}$ if $P\in\mathcal{L}(\mathscr{H})$ is a projection, while $Q\in\mathcal{L}(\mathscr{H})$ is an idempotent such that
\begin{align}\label{conditions 1 for qpp}&PQ^*P=PQP, \quad PQ^*(I-P)=-PQ(I-P),\\
\label{conditions 2 for qpp}&(I-P)Q^*(I-P)=(I-P)Q(I-P).
\end{align}
\end{definition}

Two useful characterizations of the quasi-projection pair can be stated as follows.

 \begin{lemma}\label{thm:four equivalences} \rm{ \cite[Theorem~2.4]{TXF02}} Suppose that $P\in\mathcal{L}(\mathscr{H})$ is a projection and $Q\in\mathcal{L}(\mathscr{H})$ is an idempotent.  If one element in
  \begin{equation*}\label{Sigma eight elements}\Sigma:=\big\{(A,B): A\in \{P, I-P\}, B\in \{Q,Q^*, I-Q,I-Q^*\}\big\}\end{equation*}
  is a  quasi-projection pair, then all the remaining elements in $\Sigma$ are  quasi-projection pairs.
 \end{lemma}

\begin{lemma}\label{thm:short description of qpp}\rm{ \cite[Theorem~2.5]{TXF02}} Suppose that $P\in \mathcal{L}(\mathscr{H})$ is a projection and $Q\in\mathcal{L}(\mathscr{H})$ is an idempotent. Then the following statements are equivalent:
\begin{itemize}
  \item [\rm{(i)}] $(P,Q)$ is a quasi-projection pair;
  \item [\rm{(ii)}]$Q^*=(2P-I)Q(2P-I)$;
   \item [\rm{(iii)}]$|Q^*|=(2P-I)|Q|(2P-I)$.
 \end{itemize}
\end{lemma}

It is remarkable that for every idempotent $Q$, there always exists a projection $m(Q)$, call the matched projection of $Q$, such that $\big( m(Q),Q\big)$ is a quasi-projection pair \cite[Theorem~3.1]{TXF02}. The definition of $m(Q)$ reads as follows.

\begin{definition}\cite[Definition~3.2 and Theorem~3.4]{TXF02} For each idempotent $Q\in \mathcal{L}(\mathscr{H})$, let
\begin{equation*}\label{equ:final exp of P}m(Q)=\frac12\big(|Q^*|+Q^*\big)|Q^*|^\dag\big(|Q^*|+I\big)^{-1}\big(|Q^*|+Q\big),
\end{equation*}
where $|Q^*|^\dag$ denotes the Moore-Penrose inverse of $|Q^*|$.
\end{definition}

\begin{remark}For an arbitrary idempotent  $Q\in\mathcal{L}(\mathscr{H})$, it is observed in \cite[(3.10)]{TXF02} that
$$|Q^*|^\dag=\left(P_{\mathcal{R}(Q)}P_{\mathcal{R}(Q^*)}P_{\mathcal{R}(Q)}\right)^\frac12.$$
\end{remark}

In the rest of this paper, the notation $m(Q)$ is reserved to denote the matched projection of an idempotent $Q$, and $\big(m(Q),Q\big)$ is called the matched pair of $Q$.

We end this section by recalling two basic properties of the matched projection.
\begin{lemma}\label{lem:range of matched projection}\rm{\cite[Remark~3.10]{TXF02}} If $Q\in \mathcal{L}(\mathscr{H})$ is  an idempotent, then
\begin{equation*} \mathcal{R}\big[m(Q)\big]=\mathcal{R}(|Q^*|+Q^*)=\mathcal{R}(|Q|+Q).\end{equation*}
\end{lemma}

\begin{lemma}\label{lem:property of matched projection}\rm{\cite[Theorems~3.7 and 3.14]{TXF02}} If $Q\in \mathcal{L}(\mathscr{H})$ is an idempotent, then
\begin{equation*}m(Q^*)=m(Q),\quad I-m(Q)=m(I-Q).\end{equation*}
\end{lemma}

\section{The semi-harmonious quasi-projection pairs}\label{sec:semi-harmony}

In this section, we focus mainly on the  study of the semi-harmonious quasi-projection pairs.
\subsection{Some notations associated with the quasi-projection pairs and two examples of quasi-projection pairs}
It is helpful to specify some frequently used symbols for those operators and submodules associated with two idempotents. As in \cite{TXF04}, on the whole we will introduce six operators and six submodules as follows.
\begin{definition}\label{defn of 6 operators and modules} Given idempotents $P,Q\in \mathcal{L}(\mathscr{H})$,  let
\begin{align}\label{definition of T1 and T2}&T_1=P(I-Q),  \quad    T_2=(I-P)Q,\\
\label{equ:new defn of T3 and T4}&T_3=PQ(I-P),\quad T_4=(I-P)QP,\\
&\label{equ:defn of wideT1}\widetilde{ T_1}=T_1(2P-I),\quad \widetilde{ T_2}=-T_2(2P-I),
 \end{align}
and let
\begin{align}\label{eqn:defn of H1 and H4}&\mathscr{H}_1=\mathcal{R}(P)\cap\mathcal{R}(Q), \quad \mathscr{H}_4=\mathcal{N}(P)\cap\mathcal{N}(Q),\\
 \label{eqn:defn of H2 and H3}& \mathscr{H}_2=\mathcal{R}(P)\cap\mathcal{N}(Q),  \quad \mathscr{H}_3=\mathcal{N}(P)\cap\mathcal{R}(Q), \\
 \label{eqn:defn of H5 and H6}&\mathscr{H}_5=\overline{\mathcal{R}(T_3)}, \quad  \mathscr{H}_6=\overline{\mathcal{R}(T_4)}.
\end{align}
\end{definition}

With the above notation, we establish three lemmas concerning quasi-projection pair.

\begin{lemma}\label{lem:expressions of adjoint operators}Given a quasi-projection  pair $(P,Q)$  on $\mathscr{H}$, let $T_i$ $(1\le i\le 4)$ and $\widetilde{T_i}$ $(1\le i\le 2)$ be defined by
\eqref{definition of T1 and T2}--\eqref{equ:defn of wideT1}. Then
\begin{align}&\label{exp of T1 wide T1 star}T_1^*=(2P-I)(I-Q)P,\quad \widetilde{T_1}^*=(I-Q)P ,\\
&\label{exp of T2 wide T2 star}T_2^*=-(2P-I)Q(I-P),\quad \widetilde{T_2}^*=Q(I-P), \\
&\label{exp of T3 T4 star}T_3^*=-T_4,\quad T_4^*=-T_3.
\end{align}
\end{lemma}
\begin{proof} Since $P$ ia a projection, we know that $2P-I$ is a self-adjoint unitary (symmetry). Hence, by Lemma~\ref{thm:short description of qpp}(ii) we have
\begin{align*}
   & T_2^*=(2P-I)Q(2P-I)\cdot (I-P)=-(2P-I)Q(I-P), \\
   & \widetilde{T_2}^*=-(2P-I)T_2^*=Q(I-P),\\
   & T_4^*=P\cdot (2P-I)Q(2P-I)\cdot (I-P)=-T_3.
\end{align*}
This completes the verification of \eqref{exp of T2 wide T2 star} and the second equality in \eqref{exp of T3 T4 star}.
By  Lemma~\ref{thm:four equivalences}, the remaining equalities  follow through substitution of  $(P,Q)$ by $(I-P,I-Q)$.
\end{proof}

\begin{lemma}\label{lem:the addition is self-adjoint} For every quasi-projection pair  $(P,Q)$ on $\mathscr{H}$,
  $T_1+T_2$ is self-adjoint, where $T_1$ and $T_2$ are defined by \eqref{definition of T1 and T2}.
\end{lemma}
\begin{proof} From \eqref{exp of T1 wide T1 star} and \eqref{exp of T2 wide T2 star}, we have
\begin{align*}
  (T_1+T_2)^* & =(2P-I)\big[(I-Q)P-Q(I-P)\big]=P-2PQ+Q\\
   &=P(I-Q)+(I-P)Q=T_1+T_2. \qedhere
\end{align*}
\end{proof}

 \begin{lemma}\label{lem:another form of H1 and H4} Let  $(P,Q)$  be a quasi-projection pair  on $\mathscr{H}$. Then
 \begin{align}
&\label{alternative exps of H1 and H4}\mathscr{H}_1=\mathcal{R}(P)\cap \mathcal{R}(Q^*),\quad \mathscr{H}_4=\mathcal{N}(P)\cap \mathcal{N}(Q^*),\\
&\nonumber \mathscr{H}_2=\mathcal{R}(P)\cap \mathcal{N}(Q^*),\quad \mathscr{H}_3=\mathcal{N}(P)\cap \mathcal{R}(Q^*),
\end{align}
where $\mathscr{H}_i$ $(1\le i\le 4)$ are defined by  \eqref{eqn:defn of H1 and H4}--\eqref{eqn:defn of H2 and H3}.
\end{lemma}
\begin{proof} Let $x\in \mathcal{R}(P)$. It is obvious that $(2P-I)x=x$. If furthermore $x\in \mathcal{R}(Q)$, then $Qx=x$, so by Lemma~\ref{thm:short description of qpp}(ii) we have
$$Q^*x=(2P-I)Q(2P-I)x=x.$$
It follows that $\mathscr{H}_1\subseteq \mathcal{R}(P)\cap \mathcal{R}(Q^*)$. Since $2P-I$ is a symmetry, similar reasoning gives $\mathcal{R}(P)\cap \mathcal{R}(Q^*)\subseteq \mathscr{H}_1$. This completes the verification of the first equation in
\eqref{alternative exps of H1 and H4}.

Replacing $(P,Q)$ with
$(I-P,I-Q)$, $(P,I-Q)$ and $(I-P,Q)$ respectively, the remaining three equations can be derived immediately.
\end{proof}

It is worthwhile to find out some examples of the quasi-projection pairs. For this, two examples are newly provided as follows.

\begin{example}
 Recall that an operator $J\in\mathcal{L}(\mathscr{H})$ is referred to be  a symmetry if $J^2=I$ and $J=J^*$.
 As in the Hilbert space case \cite{Ando}, we call such a pair $(\mathscr{H},J)$ (or simply $\mathscr{H}$) as a Krein space.
 An operator $Q\in\mathcal{L}(\mathscr{H})$ is said to be a weighted projection (on the indefinite inner-product space induced by $J$) if $Q^2=Q$ and $(JQ)^*=JQ$.

 Let $J_+$ be the positive part of $J$ given by $J_+=\frac12(I+J)$, and $Q\in\mathcal{L}(\mathscr{H})$ be an arbitrary idempotent. It is clear that $J_+$ is a projection, and by Lemma~\ref{thm:short description of qpp}  $(J_+,Q)$ is a quasi-projection pair if and only if $$Q^*=(2J_+-I)Q(2J_+-I),$$ or equivalently $Q^*=JQJ$, which is exactly the case that $(JQ)^*=JQ$, since $J$ is a symmetry.  This shows that $(J_+,Q)$ is a quasi-projection pair if and only if $Q$ is a weighted projection. Let $J_-=I-J_+$, which is called the negative part of $J$. From Lemma~\ref{thm:four equivalences}, we see that for every idempotent $Q\in\mathcal{L}(\mathscr{H})$, $(J_+,Q)$ is a quasi-projection pair if and only $(J_-,Q)$ is.
 \end{example}

Given a projection $P$, there might exist many idempotents $Q$ that make $(P,Q)$ to be  a quasi-projection pair.
This will be shown in our next example.

\begin{example}\label{ex:typical example of qpp}
Suppose that $A\in\mathcal{L}(\mathscr{H})_{\mbox{sa}}$ satisfying $A^2-A\geq 0$. Let
\begin{equation}\label{defn of operator ell A}\ell(A)=(A^2-A)^{\frac{1}{2}},\end{equation}
and let $P,Q\in\mathcal{L}(\mathscr{H}\oplus \mathscr{H})$ be given by
\begin{align*} P=\left(
                  \begin{array}{cc}
                    I & 0 \\
                    0 & 0 \\
                  \end{array}
                \right),\quad Q=\left(
                                  \begin{array}{cc}
                                    A & -\ell(A) \\
                                    \ell(A) & I-A \\
                                  \end{array}
                                \right).
 \end{align*}
It is a routine matter to show that $Q$ is an idempotent,  and both of \eqref{conditions 1 for qpp} and \eqref{conditions 2 for qpp} are satisfied. So, $(P,Q)$ defined as above is a quasi-projection pair.
\end{example}

\subsection{Equivalent conditions for the semi-harmonious quasi-projection pair}\label{subsec:equivalent characterizations for semi}

\begin{definition}\label{defn:semi-harmonious pair} {\rm \cite[Definition~2.2]{TXF04}} Let $(P,Q)$  be a quasi-projection pair on $\mathscr{H}$.  If $T_1$ and $T_2$ defined by \eqref{definition of T1 and T2} both have the polar decompositions, then $(P,Q)$ is said to be semi-harmonious.
\end{definition}
Suppose that $(P,Q)$ is a quasi-projection pair  on $\mathscr{H}$. As mentioned early, we will prove in Theorem~\ref{thm:simplified condition for semi-qpp} that  the above demanding for $(P,Q)$ being semi-harmonious
can in fact be weakened as  one of $T_1$ and $T_2$ has the polar decomposition.  To get such a characterization, we need two useful lemmas as follows.

\begin{lemma}\label{calculate of kernal}For every idempotents $P,Q\in\mathcal{L}(\mathscr{H})$,  the following statements are valid:
\begin{enumerate}
\item[{\rm (i)}]$\mathcal{N}\left[(I-Q)(I-P)\right]=\mathcal{R}(P)+\mathscr{H}_3$;
\item[{\rm  (ii)}]$\mathcal{N}(T_4)=\mathcal{N}(P)+\mathscr{H}_1+\mathscr{H}_2$,
\end{enumerate}
where $\mathscr{H}_1, \mathscr{H}_2,\mathscr{H}_3$ and $T_4$ are defined by \eqref{eqn:defn of H1 and H4}, \eqref{eqn:defn of H2 and H3} and \eqref{equ:new defn of T3 and T4}, respectively.
\end{lemma}
\begin{proof}A proof given in \cite[Lemma~2.1]{FXY} for two projections is also valid in the case that $P,Q\in\mathcal{L}(\mathscr{H})$ are idempotents.
\end{proof}

\begin{lemma}\label{charac of orth-PQ}
For every quasi-projection pair $(P,Q)$ on $\mathscr{H}$, let  $T_1$, $T_2$, $\widetilde{T_1}$, $\widetilde{T_2}$ and $\mathscr{H}_1$ be defined by \eqref{definition of T1 and T2}, \eqref{equ:defn of wideT1} and \eqref{eqn:defn of H1 and H4},  respectively. Then
the following statements are equivalent:
\begin{enumerate}
  \item[{\rm (i)}] $\overline{\mathcal{R}(\widetilde{T_1})}$ is orthogonally complemented in $\mathscr{H}$;
  \item[{\rm (ii)}] $\overline{\mathcal{R}(T_1)}$ is orthogonally complemented in $\mathscr{H}$;
   \item[{\rm (iii)}] $\overline{\mathcal{R}(\widetilde{T_2}^*)}$ is orthogonally complemented in $\mathscr{H}$;
  \item[{\rm (iv)}]  $\overline{\mathcal{R}(T_2^*)}$ is orthogonally complemented in $\mathscr{H}$;
  \item[{\rm (v)}]$\overline{\mathcal{R}\big[I-P+(I-Q)(I-Q)^*\big]}$ is orthogonally complemented in $\mathscr{H}$.
  \end{enumerate}
In each case, $\mathscr{H}_1$ is orthogonally complemented in $\mathscr{H}$ such that
\begin{align}\label{expression of P 1}
P=P_{\overline{\mathcal{R}(T_1)}}+P_{\mathscr{H}_1},\quad
P_{\mathcal{R}(Q)}=P_{\overline{\mathcal{R}(\widetilde{T_2}^*)}}+ P_{\mathscr{H}_1}.
\end{align}
\end{lemma}
\begin{proof}  By the first equation in \eqref{alternative exps of H1 and H4} and the second equation in \eqref{exp of T2 wide T2 star}, we have
$$\overline{\mathcal{R}(T_1)}\perp \mathscr{H}_1\quad\mbox{and}\quad \overline{\mathcal{R}(\widetilde{T_2}^*)}\perp \mathscr{H}_1.$$

(i)$\Longleftrightarrow$(ii) and (iii)$\Longleftrightarrow$(iv). From \eqref{equ:defn of wideT1} and the invertibility of $2P-I$, it is obvious that  $\mathcal{R}(\widetilde{T_1})=\mathcal{R}(T_1)$, and thus (i)$\Longleftrightarrow$(ii) follows. From \eqref{exp of T2 wide T2 star} we have
$T_2^*=(I-2P)\widetilde{T_2}^*$, which clearly leads to (iii)$\Longleftrightarrow$(iv), since  $I-2P$ is a unitary.

(v)$\Longrightarrow$(ii) and (v)$\Longrightarrow$(iii). By assumption we have $\mathscr{H}=\overline{\mathcal{R}(S)}+ \mathcal{N}(S)$, where $S$ is the positive operator given by
\begin{equation}\label{this S}S=I-P+(I-Q)(I-Q)^*.\end{equation} Clearly,
$\mathcal{N}(S)=\mathcal{R}(P)\cap \mathcal{R}(Q^*)$, which is equal to $\mathscr{H}_1$ (see \eqref{alternative exps of H1 and H4}).
Hence, $\mathscr{H}_1$ is orthogonally complemented in $\mathscr{H}$ such that
\begin{equation}\label{decomposition of H wrt S}\mathscr{H}=\overline{\mathcal{R}(S)}+ \mathscr{H}_1.\end{equation}
So each $x\in \mathscr{H}$ can be decomposed as $x=x_{1}+x_{2}$ for some
$ x_1\in \overline{\mathcal{R}(S)}$  and  $x_{2}\in \mathscr{H}_1$, which gives
\begin{equation}\label{decomposition of range of I-P}
Px=Px_1+x_2.
\end{equation}
Note that
    $$P\overline{\mathcal{R}(S)}\subseteq \overline{\mathcal{R}(PS)}=\overline{\mathcal{R}\big[(T_1(I-Q)^*\big]}\subseteq
\overline{\mathcal{R}(T_1)},$$
 so by \eqref{decomposition of range of I-P} we have
\begin{equation*}
\mathcal{R}(P)\subseteq \overline{\mathcal{R}(T_1)}+ \mathscr{H}_1\subseteq \mathcal{R}(P),
\end{equation*}
which leads to the orthogonal decomposition
\begin{equation}\label{equality of range of I-P}
\mathcal{R}(P)= \overline{\mathcal{R}(T_1)}+ \mathscr{H}_1.
\end{equation}
 Consequently,
 we arrive at
 \begin{equation*}\label{decom of (I-P)(I-Q)}
 \mathscr{H}=\overline{\mathcal{R}(T_1)}+ \overline{\mathcal{R}(T_1)}^\bot,\end{equation*}
 where
 \begin{equation*}\label{orth of (I-P)(I-Q)}
 \overline{\mathcal{R}(T_1)}^\bot=\mathcal{R}(P)^\perp+ \mathscr{H}_1.\end{equation*}
 Utilizing \eqref{equality of range of I-P} yields the first equation in \eqref{expression of P 1}. Similarly,  from \eqref{this S}, \eqref{decomposition of H wrt S} and the second equation in \eqref{exp of T2 wide T2 star}, it can be deduced that
$$\mathcal{R}(Q)=\overline{\mathcal{R}(\widetilde{T_2}^*)}+\mathscr{H}_1.$$ So, $\overline{\mathcal{R}(\widetilde{T_2}^*)}$ is orthogonally
complemented in $\mathscr{H}$ such that the second equation in \eqref{expression of P 1} is satisfied.

(iv)$\Longrightarrow$(ii). By assumption, we have
\begin{equation}\label{this range P}P\mathscr{H}=P\overline{\mathcal{R}(T_2^*)}+P\mathcal{N}(T_2).\end{equation}
A direct application of  Lemma~\ref{calculate of kernal}(i) gives
$\mathcal{N}(T_2)=\mathcal{R}(I-Q)+\mathscr{H}_1,$ which
yields
\begin{equation}\label{sth wrt P}P\mathcal{N}(T_2)=\mathcal{R}(T_1)+\mathscr{H}_1.\end{equation}  In addition, from the first equation in \eqref{exp of T2 wide T2 star} it is easily seen that
$$PT_2^*=T_1(I-P).$$
Hence, the orthogonal decomposition \eqref{equality of range of I-P} can be derived by combining the above equation with \eqref{this range P} and \eqref{sth wrt P}.
Therefore, $\overline{\mathcal{R}(T_1)}$ is orthogonally complemented in $\mathscr{H}$.

 (ii)$\Longrightarrow$(v). By assumption, we have $\mathscr{H}=\overline{\mathcal{R}(T_1)}+\mathcal{N}(T_1^*)$.
 From the first equation in \eqref{exp of T1 wide T1 star}  together with the invertibility of $2P-I$ and Lemma~\ref{calculate of kernal}(i), it can be deduced that
 $$\mathcal{N}(T_1^*)=\mathcal{N}\big[(I-Q)P\big]=\mathcal{R}(I-P)+\mathscr{H}_1.$$
 This shows that
  \begin{align}\label{decomposition of H-01}
 \mathscr{H}=\overline{\mathcal{R}(T_1)}+\mathcal{R}(I-P)+\mathscr{H}_1.
 \end{align}
Meanwhile, from Lemma~\ref{lem:Range closure of TT and T}  and \cite[Lemma~3.5]{FXY}  we can obtain
 $$\overline{\mathcal{R}(I-P)+\mathcal{R}(I-Q)}=\overline{\mathcal{R}(I-P)+\mathcal{R}\big[(I-Q)(I-Q)^*\big]}
=\overline{\mathcal{R}(S)}, $$ where $S$ is defined by \eqref{this S}. It follows that
\begin{equation*}\mathcal{R}(I-P)\subseteq \overline{\mathcal{R}(S)}\quad\mbox{and}\quad \mathcal{R}(I-Q)\subseteq \overline{\mathcal{R}(S)},\end{equation*}
which mean that  \begin{equation}\label{equ:two subsets left}\mathcal{R}(T_1)\subseteq \mathcal{R}(I-P)+\mathcal{R}(I-Q)\subseteq \overline{\mathcal{R}(S)},\end{equation}
since $T_1$ can be expressed alternatively as  $$T_1=-(I-P)(I-Q)+I-Q.$$
Combining \eqref{decomposition of H-01} and \eqref{equ:two subsets left} yields the orthogonal decomposition
\eqref{decomposition of H wrt S}.
\end{proof}

Now, we are in the position to provide some equivalent characterizations of the semi-harmonious quasi-projection pair.
\begin{theorem}\label{thm:simplified condition for semi-qpp}For every quasi-projection pair $(P,Q)$ on $\mathscr{H}$, let  $T_1$, $T_2$, $\widetilde{T_1}$, $\widetilde{T_2}$, $\mathscr{H}_1$ and $\mathscr{H}_4$ be defined by \eqref{definition of T1 and T2}, \eqref{equ:defn of wideT1} and \eqref{eqn:defn of H1 and H4}, respectively. Then the following statements are equivalent:
\begin{enumerate}
  \item [\rm{(i)}]  $(P,Q)$ is semi-harmonious;
  \item [ \rm{(ii)}] $\overline{\mathcal{R}(T_i)} (i=1,2)$ are both orthogonally complemented in $\mathscr{H}$;
  \item [ \rm{(iii)}] $\overline{\mathcal{R}(T_i^*)} (i=1,2)$ are both orthogonally complemented in $\mathscr{H}$;
  \item [ \rm{(iv)}] $\overline{\mathcal{R}\big(\widetilde{T_i}\big)} (i=1,2) $ are both orthogonally complemented in $\mathscr{H}$;
  \item [ \rm{(v)}] $\overline{\mathcal{R}\big(\widetilde{T_i}^*\big)} (i=1,2) $ are both orthogonally complemented in $\mathscr{H}$;
  \item [ \rm{(vi)}] $\overline{\mathcal{R}(T_1)}$ and $\overline{\mathcal{R}(T_1^*)}$ are both orthogonally complemented in $\mathscr{H}$;
   \item [ \rm{(vii)}] $T_1$ has the polar decomposition;
   \item [ \rm{(viii)}] $\overline{\mathcal{R}(T_1^*)}$ and $\overline{\mathcal{R}\big(\widetilde{T_2}^*\big)}$ are both orthogonally complemented in $\mathscr{H}$;
   \item [ \rm{(ix)}] $\overline{\mathcal{R}(T_2)}$ and $\overline{\mathcal{R}(T_2^*)}$ are both orthogonally complemented in $\mathscr{H}$;
    \item [ \rm{(x)}] $T_2$ has the polar decomposition;
   \item [ \rm{(xi)}] $\overline{\mathcal{R}(T_2^*)}$ and $\overline{\mathcal{R}\big(\widetilde{T_1}^*\big)}$ are both orthogonally complemented in $\mathscr{H}$.
\end{enumerate}
 In each case, $\mathscr{H}_1$ and $ \mathscr{H}_4$ are orthogonally complemented in $\mathscr{H}$ such that  \eqref{expression of P 1} is satisfied, and
\begin{equation}\label{equ:expression of PN(Q)}
I-P=P_{\overline{\mathcal{R}(T_2)}}+P_{\mathscr{H}_4},\quad
 P_{\mathcal{N}(Q)}=P_{\overline{\mathcal{R}\big(\widetilde{T_1}^*\big)}}+ P_{\mathscr{H}_4}.
\end{equation}
\end{theorem}
\begin{proof}The conclusion follows directly from Lemmas~\ref{lem:polar decomposition of T star} and \ref{charac of orth-PQ}, together with the replacement of
$(P,Q)$ by $(I-P,I-Q)$.
\end{proof}

\begin{remark}For every quasi-projection pair $(P,Q)$, equations \eqref{expression of P 1} and \eqref{equ:expression of PN(Q)} are originally  provided in
\cite[Theorem~2.6]{TXF04} with an alternative proof.
\end{remark}

\begin{corollary}\label{cor:1 to 4} For every quasi-projection pair $(P,Q)$ on $\mathscr{H}$, if one pair in
\begin{equation}\label{four pairs-pre}(P,Q),\quad (I-P,I-Q),\quad (P,Q^*) \quad\mbox{and}\quad (I-P,I-Q^*)\end{equation}
is semi-harmonious, then the remaining three pairs are all semi-harmonious.
\end{corollary}
\begin{proof} By Lemma~\ref{thm:four equivalences}, every one in \eqref{four pairs-pre} is a quasi-projection pair.
From the definition of $T_1$ and $T_2$ together with (i)$\Longleftrightarrow$(ii) in Theorem~\ref{thm:simplified condition for semi-qpp}, it is easily seen  that
\begin{align*}&\mbox{$(P,Q)$ is semi-harmonious}\Longleftrightarrow\mbox{$(I-P,I-Q)$ is semi-harmonious,}\\
&\mbox{$(P,Q^*)$ is semi-harmonious}\Longleftrightarrow\mbox{$(I-P,I-Q^*)$ is semi-harmonious.}
\end{align*}
Here the point is,
\begin{equation}\label{ranges equal wrt star-01}
  \mathcal{R}(T_1)=\mathcal{R}\big[P(I-Q)^*\big],\quad \mathcal{R}(T_2)=\mathcal{R}\big[(I-P)Q^*\big].
\end{equation}
Indeed, by Lemma~\ref{thm:short description of qpp}(ii) we have
$$(I-P)Q^*=(I-P)(2P-I)Q(2P-I)=T_2(I-2P),$$
which leads by the invertibility of $I-2P$ to the second equation in \eqref{ranges equal wrt star-01}.
Replacing $(P,Q)$ with $(I-P,I-Q)$ yields the first equation in \eqref{ranges equal wrt star-01}.
Hence, $(P,Q)$ is semi-harmonious if and only if $(P,Q^*)$ is semi-harmonious. This completes the proof.
\end{proof}

For the matched pair, we have a specific characterization as follows.

\begin{theorem}\label{thm:semi-harmony characterization for matched pair} For every idempotent $Q\in\mathcal{L}(\mathscr{H})$, its matched pair is semi-harmonious if and only if $P_{\mathcal{R}(Q)}-Q$ has the polar decomposition.
\end{theorem}
\begin{proof} What we are concerned with is the quasi-projection pair $(P,Q)$ with $P=m(Q)$. So, by \eqref{exp of T1 wide T1 star}, \eqref{exp of T2 wide T2 star} and Lemma~\ref{lem:property of matched projection} we have
\begin{equation}\label{special T 1 2 star}T_1^*=(I-Q^*)m(Q),\quad \widetilde{T_2}^*=Q\big[I-m(Q)\big]=Q\cdot m(I-Q).\end{equation}
Since $\mathcal{R}\big[m(Q)\big]=\mathcal{R}(|Q^*|+Q^*)$ (see Lemma~\ref{lem:range of matched projection}) and $\mathcal{R}(|Q^*|)=\mathcal{R}\big[P_{\mathcal{R}(Q)}\big]$, it can be deduced by  Lemma~\ref{lem:rang characterization-1} that
\begin{align*}
  \overline{\mathcal{R}(T_1^*)} & =\overline{\mathcal{R}\big[(I-Q^*)(|Q^*|+Q^*)\big]}=\overline{\mathcal{R}\big[(I-Q^*)|Q^*|\big]} \\
   & =\overline{\mathcal{R}\big[(I-Q^*)P_{\mathcal{R}(Q)}\big]}=\overline{\mathcal{R}(P_{\mathcal{R}(Q)}-Q^*)}.
\end{align*}
Replacing $Q$ with $I-Q$ in Lemma~\ref{lem:range of matched projection} gives  $$\mathcal{R}\big[m(I-Q)\big]=\mathcal{R}\big(|I-Q|+I-Q\big).$$
The equation above together with the second equation in  \eqref{special T 1 2 star} yields
 \begin{align*}
  \overline{\mathcal{R}(\widetilde{T_2}^*)} & =\overline{\mathcal{R}\big[Q(|I-Q|+I-Q)\big]} =\overline{\mathcal{R}\big[Q|I-Q|\big]}=\overline{\mathcal{R}\big[Q(I-Q^*)\big]}\\
  &=\overline{\mathcal{R}\big[QP_{\mathcal{N}(Q^*)}\big]}=\overline{\mathcal{R}\big[Q(I-P_{\mathcal{R}(Q)})\big]}=\overline{\mathcal{R}(P_{\mathcal{R}(Q)}-Q)}.
\end{align*}
Hence, by  Lemma~\ref{lem:polar decomposition of T star} we see that $P_{\mathcal{R}(Q)}-Q$ has the polar decomposition if and only if
$ \overline{\mathcal{R}(T_1^*)}$ and $\overline{\mathcal{R}(\widetilde{T_2}^*)}$ are  both orthogonally complemented in $\mathscr{H}$.
So, the conclusion can be derived from (i)$\Longleftrightarrow$(viii) in Theorem~\ref{thm:simplified condition for semi-qpp}.
\end{proof}

\subsection{A way to construct the semi-harmonious quasi-projection pairs}\label{subsec:a way}

It is worthwhile to find ways to construct the semi-harmonious quasi-projection pairs. For this, we need a lemma as follow.

\begin{lemma}\label{lem:trivial-01} Suppose that $(P,Q)$ is a quasi-projection pair on $\mathscr{H}$. Let
$T_1$ and $T_2$ be defined by \eqref{definition of T1 and T2}. Then
\begin{equation}\label{ranges equal with parameters}\mathcal{R}(\lambda_1 T_1+\lambda_2 T_2)=\mathcal{R}(T_1)+\mathcal{R}(T_2),\quad \forall \lambda_1,\lambda_2\in \mathbb{C}\setminus\{0\}.
\end{equation}
Furthermore, the following statements are equivalent:
\begin{enumerate}
\item[{\rm (i)}] $(P,Q)$ is semi-harmonious;
\item[{\rm (ii)}]$\forall\lambda_1,\lambda_2\in\mathbb{C}\setminus\{0\}$, $\overline{\mathcal{R}(\lambda_1 T_1+\lambda_2 T_2)}$ is orthogonally complemented in $\mathscr{H}$;
 \item[{\rm (iii)}]$\exists \lambda_1,\lambda_2\in\mathbb{C}\setminus\{0\}$, $\overline{\mathcal{R}(\lambda_1 T_1+\lambda_2 T_2)}$ is orthogonally complemented in $\mathscr{H}$.
 \end{enumerate}
\end{lemma}
\begin{proof} Suppose that $\lambda_1\ne 0$ and $\lambda_2\ne 0$. Let
\begin{equation*}S=\lambda_1 T_1+\lambda_2 T_2,\quad \mathscr{M}_1=\overline{\mathcal{R}(T_1)}, \quad \mathscr{M}_2=\overline{\mathcal{R}(T_2)}.
\end{equation*}
It is clear that
$$\mathcal{R}(S)\subseteq \mathcal{R}(\lambda_1 T_1)+\mathcal{R}(\lambda_2 T_2)=\mathcal{R}(T_1)+\mathcal{R}(T_2).$$
On the other hand, for every $x,y\in \mathscr{H}$ we have
$$T_1x+T_2y=S\left[\frac{1}{\lambda_1}(I-Q)x+\frac{1}{\lambda_2}Qy\right]\in \mathcal{R}(S).$$
This shows the validity of \eqref{ranges equal with parameters}.

Clearly, $\mathscr{M}_1\perp \mathscr{M}_2$, so $\mathscr{M}_1+\mathscr{M}_2$ is closed in $\mathscr{H}$. Hence, by \eqref{ranges equal with parameters}
\begin{equation}\label{closures summation}\overline{\mathcal{R}(S)}=\mathscr{M}_1+\mathscr{M}_2.\end{equation}
Therefore, $\overline{\mathcal{R}(S)}$ is invariant with respect to the parameters $\lambda_1,\lambda_2\in \mathbb{C}\setminus\{0\}$, and thus (ii)$\Longleftrightarrow$(iii) follows.

Now, assume that $\overline{\mathcal{R}(S)}$ is orthogonally complemented in $\mathscr{H}$. Then by \eqref{closures summation}
 $\mathscr{H}=\mathscr{M}_1+\mathscr{M}_1^\perp$, where $\mathscr{M}_1^\perp=\mathscr{M}_2+\overline{\mathcal{R}(S)}^\perp$.
Hence, $\mathscr{M}_1$ is orthogonally complemented in $\mathscr{H}$. The same is true for $\mathscr{M}_2$. Conversely, suppose that
$\mathscr{M}_1$ and $\mathscr{M}_2$ are both orthogonally complemented in $\mathscr{H}$. Then
$P_{\mathscr{M}_1}+P_{\mathscr{M}_2}$ is a projection whose range is equal to $\overline{\mathcal{R}(S)}$ (see \eqref{closures summation}). Therefore,
$\overline{\mathcal{R}(S)}$ is orthogonally complemented in $\mathscr{H}$. This shows that
$\overline{\mathcal{R}(S)}$ is orthogonally complemented in $\mathscr{H}$ if and only if $\mathscr{M}_1$ and $\mathscr{M}_2$
are both orthogonally complemented in $\mathscr{H}$. Hence, (i)$\Longleftrightarrow$(ii) in this lemma follows immediately from (i)$\Longleftrightarrow$(ii)
in Theorem~\ref{thm:simplified condition for semi-qpp}.
\end{proof}

Based on Lemma~\ref{lem:trivial-01}, we are now able to construct the semi-harmonious quasi-projection pairs as follows.
\begin{theorem}\label{thm:invariant submodule M}For every quasi-projection pair $(P,Q)$ on $\mathscr{H}$, let
\begin{equation}\label{equ:defn of two orthogonal parts}\mathscr{M}=\overline{\mathcal{R}(T_1)}+\overline{\mathcal{R}(T_2)},
\end{equation}
 where $T_1$ and $T_2$ are defined  by \eqref{definition of T1 and T2}.
Then   $(P|_\mathscr{M},Q|_\mathscr{M})$ is a semi-harmonious quasi-projection pair on $\mathscr{M}$.
\end{theorem}
\begin{proof}Let $\lambda_1=1$ and $\lambda_2=\pm 1$ in \eqref{ranges equal with parameters}. From \eqref{closures summation} and \eqref{definition of T1 and T2},  we obtain
\begin{equation}\label{new look at M}\mathscr{M}=\overline{\mathcal{R}(T_1+T_2)}=\overline{\mathcal{R}(T_1-T_2)}=\overline{\mathcal{R}(P-Q)},\end{equation} so $\mathscr{M}$ is closed in $\mathscr{H}$, hence it  is a Hilbert $\mathfrak{A}$-module.
Evidently,
$$\mathcal{R}(P|_\mathscr{M})=\overline{\mathcal{R}(T_1)}\subseteq \mathscr{M},$$ so $P|_\mathscr{M}$ is a projection on $\mathscr{M}$.
By \eqref{definition of T1 and T2} and Lemma~\ref{thm:short description of qpp}(ii),  we have
\begin{align*}Q(P-Q)=&PQ(P-Q)+(I-P)Q(P-Q)\\
=&P\big[I-(I-Q)\big](P-Q)+T_2(P-Q) \\
=&T_1-T_1(P-Q)+T_2(P-Q),\\
Q^*(P-Q)=&(2P-I)Q(2P-I)\cdot (P-Q) \\
=& PQ(2P-I) (P-Q)-(I-P)Q(2P-I) (P-Q) \\
=& P[I-(I-Q)](2P-I)(P-Q)-T_2(2P-I) (P-Q)\\
=& P(P-Q)-T_1(2P-I) (P-Q)-T_2(2P-I) (P-Q)\\
=&T_1-T_1(2P-I) (P-Q)-T_2(2P-I) (P-Q).
\end{align*}
From the above decompositions of  $Q(P-Q)$ and $Q^*(P-Q)$, together with \eqref{equ:defn of two orthogonal parts} and  \eqref{new look at M}, we conclude that $Q|_\mathscr{M}\in \mathcal{L}(\mathscr{M})$ such that
$(Q|_\mathscr{M})^*=Q^*|_\mathscr{M}$. Therefore, $(P|_\mathscr{M},Q|_\mathscr{M})$ is a quasi-projection pair on $\mathscr{M}$ satisfying
$$P|_\mathscr{M}(I_\mathscr{M}-Q|_\mathscr{M})=T_1|_\mathscr{M},\quad (I_\mathscr{M}-P|_\mathscr{M})Q|_\mathscr{M}=T_2|_\mathscr{M}.$$

From \eqref{equ:defn of two orthogonal parts} and the first equation in \eqref{new look at M}, it is clear that
\begin{equation*}\overline{\mathcal{R}(T_1|_\mathscr{M}+T_2|_\mathscr{M})}=\overline{\mathcal{R}\big[(T_1+T_2)|_\mathscr{M}\big]}=\overline{\mathcal{R}\big[(T_1+T_2)^2\big]}.\end{equation*}
Furthermore, $T_1+T_2$ is self-adjoint (see Lemma~\ref{lem:the addition is self-adjoint}), so from Lemma~\ref{lem:Range closure of TT and T} we have
\begin{align*}\overline{\mathcal{R}\big[(T_1+T_2)^2\big]}=\overline{\mathcal{R}\big[(T_1+T_2)(T_1+T_2)^*\big]}=\overline{\mathcal{R}(T_1+T_2)}=\mathscr{M}.
\end{align*}
This shows that $\overline{\mathcal{R}(T_1|_\mathscr{M}+T_2|_\mathscr{M})}=\mathscr{M}$, which obviously ensures the orthogonal complementarity of $\overline{\mathcal{R}(T_1|_\mathscr{M}+T_2|_\mathscr{M})}$ in $\mathscr{M}$.
The desired conclusion is immediate from (i)$\Longleftrightarrow$(iii) in Lemma~\ref{lem:trivial-01}.
\end{proof}

\begin{remark}\label{rem:prepared for the harmony}Given a quasi-projection pair $(P,Q)$ on $\mathscr{H}$, let $T_i\,(1\le i\le 4)$ and $\mathscr{M}$ be defined by
 \eqref{definition of T1 and T2}, \eqref{equ:new defn of T3 and T4} and \eqref{equ:defn of two orthogonal parts}, respectively. It is obvious that $\overline{\mathcal{R}(T_3|_\mathscr{M})}\subseteq \overline{\mathcal{R}(T_3)}$. On the other hand, by \eqref{exp of T3 T4 star} we have
\begin{align*}T_3T_3^*=-T_3T_4=-PQ(I-P)QP,
 \end{align*}
which is combined with Lemma~\ref{lem:Range closure of TT and T} and \eqref{equ:defn of two orthogonal parts} to get
 \begin{align*}\overline{\mathcal{R}(T_3)}=&\overline{\mathcal{R}(T_3T_3^*)}=\overline{\mathcal{R}\big[(PQ(I-P)QP\big]}\subseteq \overline{\mathcal{R}\big[(PQ(I-P)Q\big]}\\
 =&\overline{\mathcal{R}(T_3T_2)}=\overline{\mathcal{R}(T_3|_\mathscr{M}T_2)}\subseteq \overline{\mathcal{R}(T_3|_\mathscr{M})}.
 \end{align*}
 This shows that $\overline{\mathcal{R}(T_3|_\mathscr{M})}=\overline{\mathcal{R}(T_3)}$. Similarly, we have $\overline{\mathcal{R}(T_4|_\mathscr{M})}=\overline{\mathcal{R}(T_4)}$.
\end{remark}

\subsection{Two types of non-semi harmonious quasi-projection pairs}\label{subsec:non semi-harmonious}

In   this subsection, we provide two types of quasi-projection pairs, both of which are not semi-harmonious.

\begin{example}\label{ex:not semi-harmonious-11} Suppose that $\mathscr{H}_1$ and $\mathscr{H}_2$ are  Hilbert $C^*$-modules over a $C^*$-algebra, and
$A\in\mathcal{L}(\mathscr{H}_2,\mathscr{H}_1)$ has no polar decomposition (the reader is referred to \cite[Example~3.15]{LLX} for such an operator).
Let $Q\in\mathcal{L}(\mathscr{H}_1\oplus \mathscr{H}_2)$ be the idempotent given by
$$Q=\left(
                                                       \begin{array}{cc}
                                                         I_{\mathscr{H}_1} & A \\
                                                         0 & 0 \\
                                                       \end{array}
                                                     \right).$$
It is obvious that
$$P_{\mathcal{R}(Q)}=\left(
                       \begin{array}{cc}
                         I_{\mathscr{H}_1} & 0 \\
                         0 & 0 \\
                       \end{array}
                     \right),\quad  P_{\mathcal{R}(Q)}-Q=\left(
                                                           \begin{array}{cc}
                                                             0 & -A \\
                                                             0 & 0 \\
                                                           \end{array}
                                                         \right).$$
In virtue of Lemma~\ref{lem:polar decomposition of T star}, we see that $P_{\mathcal{R}(Q)}-Q$ has no polar decomposition. So, by Theorem~\ref{thm:semi-harmony characterization for matched pair} $\big(m(Q),Q\big)$ is not semi-harmonious.
\end{example}

To construct another type of non semi-harmonious quasi-projection pairs,  we focus our attention on the unital $C^*$-algebras. It is well-known that every $C^*$-algebra is a Hilbert module over itself. For the sake of completeness, we give a brief description. Suppose that $\mathfrak{B}$  is a $C^*$-algebra. Let
$$\langle x,y\rangle=x^*y,\quad \mbox{for $x,y\in \mathfrak{B}$}.$$
With the inner-product given as above,
 $\mathfrak{B}$ is a Hilbert $\mathfrak{B}$-module,  and $\mathfrak{B}$ can be embedded into $\mathcal{L}(\mathfrak{B})$
via $b\to L_b$ \cite[Section~3]{LMX}, where
$L_b$ is defined by $L_b(x)=bx$ for $x\in \mathfrak{B}$. If  $\mathfrak{B}$ has a unit, then it can be easily derived from \eqref{equ:keep module operator} that $\mathcal{L}(\mathfrak{B})=\{L_b:b\in \mathfrak{B}\}$. So, for every unital $C^*$-algebra $\mathfrak{B}$, we can identify $\mathfrak{B}$ with $\mathcal{L}(\mathfrak{B})$. To ease notation, we use $\mathcal{R}(b)$ and $\mathcal{N}(b)$ instead of $\mathcal{R}(L_b)$ and $\mathcal{N}(L_b)$, respectively.
Given a unital $C^*$-algebra $\mathfrak{B}$, let $\mathfrak{A}=M_2(\mathfrak{B})$, which   is a unital $C^*$-algebra.
 Hence, $\mathfrak{A}$ itself becomes a Hilbert $\mathfrak{A}$-module described as above.

\begin{remark}\label{rem:for counterexample-00}Specifically, if $\mathfrak{B}$ is a unital commutative $C^*$-algebra, then up to $C^*$-isomorphic equivalence we may assume that
$\mathfrak{B}=C(\Omega)$, the set consisting of all continuous functions on a compact Hausdorff space $\Omega$. Let $\mathfrak{A}=M_2(\mathfrak{B})$. Since $M_2(\mathbb{C})$ is nuclear and $\mathfrak{A}\cong M_2(\mathbb{C})\otimes \mathfrak{B}$, the unique $C^*$-norm on  $\mathfrak{A}$ is determined by
\begin{equation*}\label{norm of elements in c c-star alg}\|x\|=\max_{t\in \Omega}\|x(t)\|\quad (x\in \mathfrak{A}),\end{equation*}
where $x_{ij}\in C(\Omega)$ $(1\le i,j\le 2)$,
$$x=\left(
       \begin{array}{cc}
         x_{11} & x_{12} \\
         x_{21} & x_{22} \\
       \end{array}
     \right),\quad x(t)=\left(
       \begin{array}{cc}
         x_{11}(t) & x_{12}(t) \\
         x_{21}(t) & x_{22}(t) \\
       \end{array}
     \right)\quad (t\in\Omega),$$
and for each $t\in\Omega$, $\|x(t)\|$
denotes the  operator norm (also known as 2-norm) on $M_2(\mathbb{C})$.
\end{remark}

\begin{theorem}\label{thm:for counterexample-01}Suppose that $\mathfrak{B}$ is a unital $C^*$-algebra and $A\in \mathfrak{B}_{\mbox{sa}}$ is such that $A^2-A\geq 0$.
Let $\mathfrak{A}=M_2(\mathfrak{B})$ and $\ell(A)\in \mathfrak{B}$ be given by
 \eqref{defn of operator ell A},
  and let $P,Q\in \mathfrak{A}$ be defined by
 \begin{align}\label{equ:nnew defn of P Q unital Cstar alg} P=\left(
                  \begin{array}{cc}
                    I_{\mathfrak{B}} & 0 \\
                    0 & 0 \\
                  \end{array}
                \right),\quad Q=\left(
                                  \begin{array}{cc}
                                    A & -\ell(A) \\
                                    \ell(A) & I_{\mathfrak{B}}-A \\
                                  \end{array}
                                \right).
 \end{align}
 Then $(P,Q)$ is semi-harmonious if and only if
  \begin{equation}\label{equ:for semi-harmony in unital alg}I_{\mathfrak{B}}\in \overline{\mathcal{R}(A-I_{\mathfrak{B}})}+\mathcal{N}(A-I_{\mathfrak{B}}),
 \end{equation}
 where
 \begin{align*}&\mathcal{R}(A-I_{\mathfrak{B}})=\big\{(A-I_{\mathfrak{B}})x:x\in \mathfrak{B}\big\},\\
 &\mathcal{N}(A-I_{\mathfrak{B}})=\big\{x\in \mathfrak{B}:(A-I_{\mathfrak{B}})x=0\big\}.
 \end{align*}
\end{theorem}
 \begin{proof} From Example~\ref{ex:typical example of qpp}, we know that $(P,Q)$ is a quasi-projection pair on $\mathfrak{A}$.
 It is clear that  the operators $T_1$ and $T_2$ defined by \eqref{definition of T1 and T2} turn out to be
 \begin{equation}\label{this new expression for T 1 2}T_1=\left(
         \begin{array}{cc}
           I_{\mathfrak{B}}-A & \ell(A) \\
           0 & 0 \\
         \end{array}
       \right),\quad T_2=\left(
         \begin{array}{cc}
           0 & 0 \\
           \ell(A) & I_{\mathfrak{B}}-A \\
         \end{array}
       \right).
 \end{equation}
 Let \begin{equation*}\mathscr{X}=\mathcal{R}(A-I_{\mathfrak{B}}),\quad \mathscr{Y}=\mathcal{N}(A-I_{\mathfrak{B}}).
 \end{equation*}
  Direct computations yield
 $$(I_{\mathfrak{B}}-A)^2+\ell^2(A)=(A-I_{\mathfrak{B}})(2A-I_{\mathfrak{B}})=(2A-I_{\mathfrak{B}})(A-I_{\mathfrak{B}}).$$
 By assumption $A^2-A\geq 0$, which implies that $2A-I_{\mathfrak{B}}$ is invertible in $\mathfrak{B}$. So, from the above equations
  we obtain
 \begin{align*}&\mathcal{R}\big[(I_{\mathfrak{B}}-A)^2+\ell^2(A)\big]=\mathscr{X},\quad \mathcal{N}\big[(I_{\mathfrak{B}}-A)^2+\ell^2(A)\big]=\mathscr{Y}.
 \end{align*}
  Therefore,
 \begin{align*}&\overline{\mathcal{R}(T_1)}=\overline{\mathcal{R}(T_1T_1^*)}
 =\left\{\left(
           \begin{array}{cc}
             x_{11} & x_{12} \\
             0 & 0 \\
           \end{array}
         \right): x_{11}, x_{12}\in \overline{\mathscr{X}}\right\},\\
  &\mathcal{N}(T_1^*)=\mathcal{N}(T_1T_1^*)=\left\{\left(
           \begin{array}{cc}
             x_{11} & x_{12} \\
             x_{21} & x_{22} \\
           \end{array}
         \right): x_{11}, x_{12}\in \mathscr{Y}; x_{21}, x_{22}\in \mathfrak{B}\right\}.
 \end{align*}
The equations above together with \eqref{equ:keep module operator} indicate that
 $$\mathfrak{A}=\overline{\mathcal{R}(T_1)}+\mathcal{N}(T_1^*)\Longleftrightarrow \mathfrak{B}=\overline{\mathscr{X}}+\mathscr{Y}\Longleftrightarrow I_{\mathfrak{B}}\in \overline{\mathscr{X}}+\mathscr{Y}.$$
This shows that $\overline{\mathcal{R}(T_1)}$ is orthogonally complemented in $\mathfrak{A}$ if and only if \eqref{equ:for semi-harmony in unital alg} is satisfied.
From \eqref{this new expression for T 1 2} and the derivations as above, it is easily seen that the same is true if $T_1$ is replaced by $T_2$. The desired conclusion follows from
(i)$\Longleftrightarrow$(ii) in Theorem~\ref{thm:simplified condition for semi-qpp}.
\end{proof}

In view of  Remark~\ref{rem:for counterexample-00} and Theorem~\ref{thm:for counterexample-01}, it is easy to construct new type of non-semi harmonious quasi-projection pairs.
We provide an example as follows.

\begin{example}\label{prop:counter example-1}
 Let $\mathfrak{B}=C(\Omega)$ with $\Omega=[0,1]$ and $\mathfrak{A}=M_2(\mathfrak{B})$, which is
regarded as the Hilbert module over itself. Let $P$ and $Q$ be defined by \eqref{equ:nnew defn of P Q unital Cstar alg} such that
the element $A$ in $\mathfrak{B}$ is given  by
$$A(t)=\sec^2(t),\quad t\in [0,1].$$  Then
$(A-I_{\mathfrak{B}})(t)=\tan^2(t)$ for any $t\in [0,1]$, whence
$$\mathcal{N}(A-I_{\mathfrak{B}})=\big\{x\in C[0,1]:\tan^2(t)x(t)\equiv 0, t\in [0,1]\big\}=\{0\},$$
where $x(t)=0$ for $t\ne 0$ comes from the observation that $\tan^2(t)\ne 0$ for all such $t$, and $x(0)=0$ is derived from the continuity of $x$ at $t=0$.

By definition, $\mathcal{R}(A-I_{\mathfrak{B}})=\{(A-I_{\mathfrak{B}})y:y\in \mathfrak{B}\}=\{\tan^2(\cdot )y:y\in \mathfrak{B}\}$.
In virtue of $\tan(0)=0$, we have $x(0)=0$ for every $x\in
\overline{\mathcal{R}(A-I_{\mathfrak{B}})}$.
Therefore, \eqref{equ:for semi-harmony in unital alg} fails to be true. Hence by Theorem~\ref{thm:for counterexample-01}, $(P,Q)$ is not semi-harmonious.
\end{example}

\section{The harmonious quasi-projection pairs}\label{sec:harmony}
The purpose of this section is to set up the general theory for the harmonious quasi-projection pairs, which will play a crucial role in  the block matrix representations for the quasi-projection pairs \cite{TXF05}.
Inspired by \cite[Definition~4.1]{LMX} and \cite[Theorem~2.4]{FXY}, we make a definition as follows.
\begin{definition}\label{defn of harmonious quasi-projection pair} A quasi-projection pair $(P,Q)$ on $\mathscr{H}$ is referred to be harmonious if
$\mathscr{H}_5$ and $\mathscr{H}_6$ defined by \eqref{eqn:defn of H5 and H6} are both orthogonally complemented in $\mathscr{H}$.
\end{definition}

\begin{remark}\label{rem:simplified condition for harmonious --0}Given a quasi-projection pair $(P,Q)$ on $\mathscr{H}$, let $T_3$ and $T_4$ be defined by \eqref{equ:new defn of T3 and T4}.
We may combine Definition~\ref{defn of harmonious quasi-projection pair} with \eqref{eqn:defn of H5 and H6}, \eqref{exp of T3 T4 star} and
Lemma~\ref{lem:polar decomposition of T star} to conclude that
\begin{align*}\mbox{$(P,Q)$ is harmonious}&\Longleftrightarrow \mbox{$T_3$ has the polar decomposition}\\
&\Longleftrightarrow \mbox{$T_4$ has the polar decomposition}.
\end{align*}
\end{remark}

As literally revealed, every harmonious quasi projection pair must be a semi-harmonious quasi projection pair. Actually, a more in-depth characterization will be provided in Theorem~\ref{thm:harmony=2 semi-harmony}. To achieve the desired characterization, we need a lemma as follows.

 \begin{lemma}\label{charac of orth-PQ(I-P)}
 For every quasi-projection pair $(P,Q)$ on $\mathscr{H}$, let $T_1$ and $\mathscr{H}_5$ be defined by \eqref{definition of T1 and T2} and \eqref{eqn:defn of H5 and H6}, respectively.
  Then the following statements are equivalent:
\begin{enumerate}
  \item[{\rm (i)}] $\mathscr{H}_5$ is orthogonally complemented in $\mathscr{H}$;
  \item[{\rm (ii)}] $\overline{\mathcal{R}(PQ)}$  and $\overline{\mathcal{R}(T_1)}$ are both
       orthogonally complemented in $\mathscr{H}$.
\end{enumerate}
In each case, $\mathscr{H}_1$ and $\mathscr{H}_2$ are both orthogonally complemented in $\mathscr{H}$ such that
\begin{align}\label{decomposition of PQ in two steps}&P_{\overline{\mathcal{R}(PQ)}}=P_{\mathscr{H}_1}+P_{\mathscr{H}_5},\quad  P=P_{\overline{\mathcal{R}(PQ)}}+P_{\mathscr{H}_2},\quad P_{\overline{\mathcal{R}(T_1)}}=P_{\mathscr{H}_2}+P_{\mathscr{H}_5}.
\end{align}
 \end{lemma}
 \begin{proof} Let $T_3, T_4, \mathscr{H}_1$ and $\mathscr{H}_2$ be defined by \eqref{equ:new defn of T3 and T4}, \eqref{eqn:defn of H1 and H4} and \eqref{eqn:defn of H2 and H3}, respectively.

  (i)$\Longrightarrow$(ii). Assume  that $\mathscr{H}=\mathscr{H}_5+\mathscr{H}_5^\perp$. From \eqref{eqn:defn of H5 and H6} and \eqref{exp of T3 T4 star}, we have
 $$\mathscr{H}_5^\perp=\overline{\mathcal{R}(T_3)}^\perp=\mathcal{N}(T_3^*)=\mathcal{N}(T_4).$$
 So, a direct use of Lemma~\ref{calculate of kernal}(ii) gives
 \begin{equation}\label{orthogonal decomposition of H into 4 parts}\mathscr{H}=\mathscr{H}_1+\mathscr{H}_2+\mathscr{H}_5+\mathcal{N}(P).\end{equation}
By the definitions of $\mathscr{H}_1,\mathscr{H}_2$ and $\mathscr{H}_5$, it is clear that
$$\mathcal{N}(P)\perp \mathscr{H}_i (i=1,2,5),\quad \mathscr{H}_5\perp \mathscr{H}_j (j=1,2).$$
Combining \eqref{alternative exps of H1 and H4} and  \eqref{eqn:defn of H2 and H3} yields $\mathscr{H}_1\perp \mathscr{H}_2$. Hence, it can be concluded by \eqref{orthogonal decomposition of H into 4 parts} that $\mathscr{H}_i(i=1,2,5)$ are all orthogonally complemented in $\mathscr{H}$ such that
\begin{equation}\label{orthogonal decomposition of P into 3 parts}P=P_{\mathscr{H}_1}+P_{\mathscr{H}_2}+P_{\mathscr{H}_5}.\end{equation}

Rewrite \eqref{orthogonal decomposition of H into 4 parts} as
$$\mathscr{H}=\mathscr{X}_1+\mathscr{X}_2=\mathscr{Y}_1+\mathscr{Y}_2,$$
where $\mathscr{X}_i, \mathscr{Y}_i (i=1,2)$ are defined by
$$\mathscr{X}_1=\mathscr{H}_2+\mathscr{H}_5,\quad \mathscr{X}_2=\mathscr{H}_1+\mathcal{N}(P),\quad \mathscr{Y}_1=\mathscr{H}_1+\mathscr{H}_5,\quad \mathscr{Y}_2=\mathscr{H}_2+\mathcal{N}(P).$$
Note that $2P-I$ is surjective and $P(2P-I)=P$, so by Lemma~\ref{thm:short description of qpp}(ii) we have $\mathcal{R}(PQ^*)=\mathcal{R}(PQ)$.
It follows  that
\begin{align*}&\mathscr{Y}_1\subseteq \overline{\mathcal{R}(PQ)}, \quad  \mathscr{Y}_2\subseteq \mathcal{N}(QP)=\overline{\mathcal{R}(PQ^*)}^\bot=\overline{\mathcal{R}(PQ)}^\bot.
  \end{align*}
This shows that
$$\mathscr{H}= \overline{\mathcal{R}(PQ)}+\overline{\mathcal{R}(PQ)}^\bot,$$
hence $\overline{\mathcal{R}(PQ)}$ is orthogonally complemented in $\mathscr{H}$. Moreover,
\begin{align*}&P_{\overline{\mathcal{R}(PQ)}}\mathscr{H}=P_{\overline{\mathcal{R}(PQ)}}(\mathscr{Y}_1+\mathscr{Y}_2)=P_{\overline{\mathcal{R}(PQ)}}\mathscr{Y}_1=\mathscr{Y}_1=(P_{\mathscr{H}_1}+P_{\mathscr{H}_5})\mathscr{H},
\end{align*}
so, the first equation in \eqref{decomposition of PQ in two steps} is satisfied. In view of \eqref{orthogonal decomposition of P into 3 parts},
the second  equation in \eqref{decomposition of PQ in two steps} follows.
Similarly,
 \begin{align*}\mathscr{X}_1\subseteq \overline{\mathcal{R}(T_1)},\quad \mathscr{X}_2\subseteq \mathcal{N}\big[(I-Q)P\big]=\overline{\mathcal{R}\big[P(I-Q)^*\big]}^\bot=\overline{\mathcal{R}(T_1)}^\bot.
  \end{align*}
Therefore, $\overline{\mathcal{R}(T_1)}$ is orthogonally complemented in $\mathscr{H}$ such that the last equation in \eqref{decomposition of PQ in two steps} is satisfied.

(ii)$\Longrightarrow$(i). From the first equation in \eqref{expression of P 1}, we have
\begin{equation}\label{decom of H in PQ(I-P)} \mathscr{H}=\overline{\mathcal{R}(T_1)}+\mathcal{R}(I-P)+\mathscr{H}_1.\end{equation}
By a replacement of $(P,Q)$ with $(P,I-Q)$, we obtain
\begin{align}\label{decom of H in PQ}\mathscr{H}=&\overline{\mathcal{R}(PQ)}+\mathcal{N}(P)+\mathscr{H}_2.
\end{align}
It follows from \eqref{decom of H in PQ(I-P)} and \eqref{equ:new defn of T3 and T4} that
\begin{align*}
\mathcal{R}(PQ)\subseteq & \overline{\mathcal{R}\left[PQT_1\right]}+\mathcal{R}(T_3)+\mathscr{H}_1\\ \nonumber
=&\overline{\mathcal{R}\left[T_3(I-Q)\right]}+\mathcal{R}(T_3)+\mathscr{H}_1 \subseteq \mathscr{H}_5+\mathscr{H}_1.
\end{align*}
Therefore,
\begin{equation*}
\overline{\mathcal{R}(PQ)}\subseteq \mathscr{H}_5+\mathscr{H}_1.
\end{equation*}
This together with \eqref{decom of H in PQ} yields  the orthogonal decomposition $\mathscr{H}=\mathscr{H}_5+ \mathscr{Y}$, where $\mathscr{Y}$ is the right side
of Lemma~\ref{calculate of kernal}(ii).
\end{proof}

\begin{theorem}\label{thm:harmony=2 semi-harmony} For every quasi-projection pair $(P,Q)$ on $\mathscr{H}$, $(P,Q)$ is harmonious if and only if
$(P,Q)$ and $(P, I-Q)$ are both semi-harmonious.
\end{theorem}
\begin{proof}By a replacement of $(P,Q)$ with $(I-P,I-Q)$, it can be deduced from Lemma~\ref{charac of orth-PQ(I-P)} that
$\mathscr{H}_6$ defined by \eqref{eqn:defn of H5 and H6} is orthogonally complemented in $\mathscr{H}$ if and only if
$\overline{\mathcal{R}\big[(I-P)(I-Q)\big]}$  and $\overline{\mathcal{R}(T_2)}$ are both
       orthogonally complemented in $\mathscr{H}$. Thus, $(P,Q)$ is harmonious if and only if
$\overline{\mathcal{R}(T_1)}$, $\overline{\mathcal{R}(T_2)}$, $\overline{\mathcal{R}(PQ)}$ and $\overline{\mathcal{R}\big[(I-P)(I-Q)\big]}$  are all orthogonally complemented in $\mathscr{H}$.
So, the desired conclusion  follows from (i)$\Longleftrightarrow$(ii) in Theorem~\ref{thm:simplified condition for semi-qpp}.
\end{proof}

As a consequence of Theorem~\ref{thm:harmony=2 semi-harmony}, we have the following symmetrical result.
\begin{corollary}\label{cor:four equivalences} For every quasi-projection pair  $(P,Q)$ on $\mathscr{H}$, let
$\Sigma$ be defined by \eqref{Sigma eight elements}. If one element in $\Sigma$
  is harmonious, then all the remaining elements in $\Sigma$ are harmonious.
 \end{corollary}
\begin{proof}We may as well assume that $(P,Q)$ is harmonious. In this case, by Theorem~\ref{thm:harmony=2 semi-harmony} $(P,Q)$ and $(P, I-Q)$ are both semi-harmonious. Repeatedly using Corollary~\ref{cor:1 to 4}  twice indicates that the four pairs in \eqref{four pairs-pre} and these in
$$(P,I-Q),\quad (I-P,Q),\quad (P,I-Q^*) \quad\mbox{and}\quad (I-P,Q^*)$$
are all semi-harmonious. Hence, the conclusion is  immediate from Theorem~\ref{thm:harmony=2 semi-harmony}.
\end{proof}

\begin{remark}\label{charac of orth-H6} Suppose that  $(P,Q)$ is a harmonious quasi-projection pair. Let
$T_1, T_2$ and $\mathscr{H}_i (1\le i\le 6)$ be defined by \eqref{definition of T1 and T2}, \eqref{eqn:defn of H1 and H4}, \eqref{eqn:defn of H2 and H3} and \eqref{eqn:defn of H5 and H6}, respectively. By Definition~\ref{defn of harmonious quasi-projection pair}, $\mathscr{H}_5$ and $\mathscr{H}_6$ are orthogonally complemented in $\mathscr{H}$. So, by Lemma~\ref{charac of orth-PQ(I-P)} $\mathscr{H}_1$ and $\mathscr{H}_2$ are orthogonally complemented in $\mathscr{H}$ such that \eqref{decomposition of PQ in two steps} is satisfied. Replacing $(P,Q)$ with $(I-P,Q)$,
we see that $\mathscr{H}_3$ and $\mathscr{H}_4$  are also orthogonally complemented in $\mathscr{H}$ such that
\begin{align*}
   & P_{\overline{\mathcal{R}(T_2)}}=P_{\mathscr{H}_3}+P_{\mathscr{H}_6},\quad I-P=P_{\overline{\mathcal{R}(T_2)}}+P_{\mathscr{H}_4}, \\
\label{decomposition of I-P in two steps}&P_{\overline{\mathcal{R}[(I-P)(I-Q)]}}=P_{\mathscr{H}_4}+P_{\mathscr{H}_6}.
\end{align*}
\end{remark}

We provide a featured result of the matched pair as follows.
\begin{theorem}\label{thm:matched pair semi harmony is harmony} For every matched pair $\big(m(Q),Q\big)$, it is harmonious if and only if it is semi-harmonious.
\end{theorem}
\begin{proof} Let $Q$ be an arbitrary idempotent on $\mathscr{H}$. By Theorem~\ref{thm:harmony=2 semi-harmony},
$\big(m(Q),Q\big)$ is harmonious if and only if both $\big(m(Q),Q\big)$ and $\big(m(Q),I-Q\big)$ are semi-harmonious. So, to get the conclusion, it is sufficient to prove that  $\big(m(Q),I-Q\big)$ is  always semi-harmonious.
Let $$T_1=m(Q)Q,\quad T_2=\big[I-m(Q)\big](I-Q).$$ From Lemma~\ref{lem:range of matched projection}, Lemma~\ref{lem:rang characterization-1} and the invertibility of $|Q^*|+I$, we have
$$ \overline{\mathcal{R}(T_1^*)}=\overline{\mathcal{R}\big[Q^*(|Q^*|+Q^*)\big]}=\overline{\mathcal{R}\big[Q^*(|Q^*|+I)\big]}=\mathcal{R}(Q^*),$$
which is orthogonally complemented in $\mathscr{H}$. In virtue of the second equation stated in Lemma~\ref{lem:property of matched projection}, it can be concluded that $\overline{\mathcal{R}(T_2^*)}=\mathcal{R}(I-Q^*)$  by a simple replacement of $Q$ with $I-Q$.
Therefore, $\overline{\mathcal{R}(T_2^*)}$ is also orthogonally complemented in $\mathscr{H}$.
Hence, the desired conclusion can be derived immediately from (i)$\Longleftrightarrow$(iii) in Theorem~\ref{thm:simplified condition for semi-qpp}.
\end{proof}

As we know, every closed linear subspace of a Hilbert space is always orthogonally complemented. So, if $(P,Q)$ is a quasi-projection pair acting on a Hilbert space, then it is harmonious. It is notable that the same is true for some well-behaved Hilbert $C^*$-modules \cite{Magajna}.

To construct new examples of the harmonious and non-harmonious quasi-projection pairs, it is helpful to study furthermore the unital $C^*$-algebras, which are regarded as Hilbert $C^*$-modules.

\begin{theorem}\label{thm:for harmony in unital alg} Suppose that $\mathfrak{B}$ is a unital $C^*$-algebra and $A\in \mathfrak{B}_{\mbox{sa}}$ is such that $A^2-A\geq 0$.
Let $\mathfrak{A}=M_2(\mathfrak{B})$ and $\ell(A)\in \mathfrak{B}$ be given by
 \eqref{defn of operator ell A},
  and let $P,Q\in \mathfrak{A}$ be defined by \eqref{equ:nnew defn of P Q unital Cstar alg}.
Then $(P,Q)$ is harmonious if and only if
 \begin{equation}\label{equ:for harmony in unital alg}I_{\mathfrak{B}}\in \overline{\mathcal{R}(A^2-A)}+\mathcal{N}(A^2-A),
 \end{equation}
 where
\begin{align*}&\mathcal{R}(A^2-A)=\big\{(A^2-A)x:x\in \mathfrak{B}\big\},\\
 &\mathcal{N}(A^2-A)=\big\{x\in \mathfrak{B}: (A^2-A)x=0\big\}.
 \end{align*}
   \end{theorem}
 \begin{proof} Let $T_3$ and $T_4$ be defined by \eqref{equ:new defn of T3 and T4}. Direct computations yield
 \begin{equation}\label{nnew expression for T3 T4}T_3=\left(
         \begin{array}{cc}
           0 & -\ell(A) \\
           0 & 0 \\
         \end{array}
       \right),\quad T_4=\left(
         \begin{array}{cc}
           0 & 0 \\
           \ell(A) & 0 \\
         \end{array}
       \right),
 \end{equation}
 in which $\ell(A)$ is given by \eqref{defn of operator ell A}.
 Put \begin{align*}\mathscr{X}=\mathcal{R}\big(\ell(A)\big),\quad \mathscr{Y}=\mathcal{N}\big(\ell(A)\big).
 \end{align*}
 Since $\ell(A)\ge 0$, by Lemmas~\ref{lem:Range closure of TT and T} and \ref{lem:Range Closure of T alpha and T} we have
 \begin{align*}&\overline{\mathscr{X}}=\overline{\mathcal{R}\big(\ell(A)\big)}=\overline{\mathcal{R}\big(\ell^2(A)\big)}=\overline{\mathcal{R}(A^2-A)},\quad \mathscr{Y}=\mathcal{N}\big(\ell^2(A)\big)=\mathcal{N}(A^2-A).
 \end{align*}
 It follows from \eqref{nnew expression for T3 T4} that
 \begin{align*}&\mathscr{H}_5=\overline{\mathcal{R}(T_3)}=\left\{\left(
           \begin{array}{cc}
             x_{11} & x_{12} \\
             0 & 0 \\
           \end{array}
         \right): x_{11}, x_{12}\in \overline{\mathscr{X}}\right\},\\
  &\mathscr{H}_5^\perp=\mathcal{N}(T_3^*)=\left\{\left(
           \begin{array}{cc}
             x_{11} & x_{12} \\
             x_{21} & x_{22} \\
           \end{array}
         \right): x_{11}, x_{12}\in \mathscr{Y}; x_{21}, x_{22}\in \mathfrak{B}\right\}.
 \end{align*}
Hence,
 $$\mathfrak{A}=\mathscr{H}_5+\mathscr{H}_5^\perp\Longleftrightarrow \mathfrak{B}=\overline{\mathscr{X}}+\mathscr{Y}\Longleftrightarrow I_{\mathfrak{B}}\in \overline{\mathscr{X}}+\mathscr{Y}.$$
This shows that $\mathscr{H}_5$ is orthogonally complemented in $\mathscr{H}$ if and only if \eqref{equ:for harmony in unital alg} is satisfied.
Similarly, it can also be concluded by \eqref{nnew expression for T3 T4} that $\mathscr{H}_6$ is orthogonally complemented in $\mathscr{H}$ if and only if \eqref{equ:for harmony in unital alg} is satisfied.
\end{proof}

\begin{example}\label{ex:change to harmony}Let $\mathfrak{B}$, $\mathfrak{A}$, $P$ and $Q$ be the same as in Example~\ref{prop:counter example-1}. As shown in Example~\ref{prop:counter example-1},
$(P,Q)$ fails to be semi-harmonious. Denote $I_{\mathfrak{A}}$ simply by $I$, and let $T_1,T_2, T_3,T_4$ and $\mathscr{M}$ be defined by
 \eqref{definition of T1 and T2}, \eqref{equ:new defn of T3 and T4} and \eqref{equ:defn of two orthogonal parts}, respectively. According to Theorem~\ref{thm:invariant submodule M},
$(P|_\mathscr{M},Q|_\mathscr{M})$ is semi-harmonious. In what follows,  we  prove that $(P|_\mathscr{M},Q|_\mathscr{M})$ is in fact a harmonious quasi-projection pair.

It should be aware that for a function $x$ in $C[0,1]$, the notation $\mathcal{R}(x)$ is usually employed for the range $\{x(t):t\in [0,1]\}$ of $x$.
To avoid such an ambiguity, it is helpful to use
$\mathcal{R}\big(x(\cdot)\big)$ to denote the set $\big\{xy:y\in C[0,1]\big\}$.
 Note that $\sec(t)\ne 0$ for every $t\in [0,1]$, so $\mathcal{R}\big[\sec(\cdot)\tan(\cdot)\big]=\mathcal{R}\big(\tan(\cdot)\big)$.
Therefore, by Lemma~\ref{lem:Range closure of TT and T} we have
\begin{equation}\label{equ:derivation of equal to X}\mathscr{X}:=\overline{\mathcal{R}\big[\sec(\cdot)\tan(\cdot)\big]}=\overline{\mathcal{R}\big(\tan^2(\cdot)\big)}.
\end{equation}

Let $C,D\in\mathcal{L}(\mathfrak{B})$ be defined by
 \begin{equation*}\label{equ:new defn of C and D-self-adjoint}C=P(I-Q)P,\quad D=(I-P)Q(I-P).
 \end{equation*}
It follows immediately  from \eqref{conditions 1 for qpp} and \eqref{conditions 2 for qpp} that
\begin{equation*}\label{equ:C and D are hermitian}C=C^*\quad\mbox{and}\quad D=D^*.
\end{equation*}
Since $C=T_1P$, we have $\overline{\mathcal{R}(C)}\subseteq \overline{\mathcal{R}(T_1)}$.
Utilizing  \eqref{definition of T1 and T2} and \eqref{exp of T1 wide T1 star} yields
\begin{equation*}\label{equ:T1 multiplied by T1 star}T_1T_1^*=P(I-Q)(2P-I)(I-Q)P=2C^2-C.
 \end{equation*}
 So, according to Lemma~\ref{lem:Range closure of TT and T} we have
 \begin{align*}\overline{\mathcal{R}(T_1)}=&\overline{\mathcal{R}(T_1T_1^*)}=\overline{\mathcal{R}(2C^2-C)}\subseteq \overline{\mathcal{R}(C)}.
 \end{align*}
 Consequently,
 \begin{equation}\label{equ:closure of range T1-three terms}\overline{\mathcal{R}(T_1)}=\overline{\mathcal{R}(2C^2-C)}=\overline{\mathcal{R}(C)}.
 \end{equation}
 Replacing $(P,Q)$ with $(I-P,I-Q)$ gives immediately
 \begin{equation}\label{equ:closure of range T2-three terms}T_2T_2^*=2D^2-D,\quad \overline{\mathcal{R}(T_2)}=\overline{\mathcal{R}(2D^2-D)}=\overline{\mathcal{R}(D)}.
 \end{equation}
 In view of \eqref{equ:defn of two orthogonal parts}, \eqref{equ:closure of range T1-three terms} and \eqref{equ:closure of range T2-three terms},  we have $\mathscr{M}=\overline{\mathscr{M}_0}$, where
 $$\mathscr{M}_0=\mathcal{R}(C)+\mathcal{R}(D).$$
A direct computation shows that for each $t\in [0,1]$,
\begin{equation*}C(t)=\left(
          \begin{array}{cc}
            -\tan^2(t) & 0 \\
            0 &  0\\
          \end{array}
        \right),\quad  D(t)=\left(
          \begin{array}{cc}
            0 & 0 \\
            0 & -\tan^2(t) \\
          \end{array}
        \right),
\end{equation*}
which means that for any $x\in \mathfrak{A}$ with $x(t)=\big(x_{ij}(t)\big)_{1\le i,j\le 2}$,
$$x\in \mathscr{M}_0\Longleftrightarrow x_{ij}\in\mathcal{R}\big(\tan^2(\cdot)\big) \quad\mbox{for $i,j=1,2$}.$$
Hence,
\begin{equation}\label{equ:characterization of x in X}x\in \mathscr{M}\Longleftrightarrow x_{ij}\in \mathscr{X}\quad\mbox{for $i,j=1,2$},
\end{equation}
where $\mathscr{X}$ is defined by \eqref{equ:derivation of equal to X}.

From Remark~\ref{rem:prepared for the harmony} and \eqref{exp of T3 T4 star}, we see that the verification of the orthogonal decomposition
$$\mathscr{M}=\overline{\mathcal{R}(T_4|_\mathscr{M})}+\mathcal{N}(T_4^*|_\mathscr{M})$$
can be simplified as
\begin{equation}\label{equ:new orth decomp of M}\mathscr{M}=\overline{\mathcal{R}(T_4)}+\mathcal{N}(T_3)\cap \mathscr{M}.
\end{equation}
For every $x\in \mathfrak{A}$ with $x(t)=\big(x_{ij}(t)\big)_{1\le i,j\le 2}$, by \eqref{nnew expression for T3 T4} it is easy to verify that
$$x\in \mathcal{N}(T_3)\Longleftrightarrow x_{21}(t)\equiv 0,\, x_{22}(t)\equiv 0.$$
So, it can be deduced from \eqref{equ:characterization of x in X} that
\begin{equation}\label{equ:for the derivation of orth decomp-01}x\in \mathcal{N}(T_3)\cap \mathscr{M}\Longleftrightarrow x_{21}(t)\equiv 0, x_{22}(t)\equiv 0; x_{11}\in \mathscr{X}, x_{12}\in \mathscr{X}.\end{equation}
Substituting $A(t)=\sec^2(t)$ into \eqref{nnew expression for T3 T4} yields
$$T_4(t)=\left(
           \begin{array}{cc}
             0 & 0 \\
             \sec(t)\tan(t) & 0 \\
           \end{array}
         \right)$$
for $t\in [0,1]$. Therefore, it can be concluded by \eqref{equ:derivation of equal to X} that
\begin{equation}\label{equ:for the derivation of orth decomp-02}x\in \overline{\mathcal{R}(T_4)}\Longleftrightarrow x_{11}(t)\equiv 0, x_{12}(t)\equiv 0; x_{21}\in \mathscr{X}, x_{22}\in \mathscr{X}.\end{equation}
Evidently, \eqref{equ:new orth decomp of M} can be derived immediately from \eqref{equ:characterization of x in X}, \eqref{equ:for the derivation of orth decomp-01} and \eqref{equ:for the derivation of orth decomp-02}.
This shows that  $\overline{\mathcal{R}(T_4|_\mathscr{M})}$ is orthogonally complemented in $\mathscr{M}$. Similar reasoning ensures the orthogonal complementarity of
  $\overline{\mathcal{R}(T_3|_\mathscr{M})}$ in $\mathscr{M}$.
\end{example}

There exist semi-harmonious quasi-projection pairs, which fail to be  harmonious. Our first example in this direction is as follows.
\begin{example}\label{ex:semi-harmony does not imply harmony-01}
Let $Q\in\mathcal{L}(\mathscr{H})$ be an idempotent such that $P_{\mathcal{R}(Q)}-Q$ has no polar decomposition. In this case,
by Theorem~\ref{thm:semi-harmony characterization for matched pair} $\big(m(Q),Q\big)$ is not semi-harmonious,
hence by Theorem~\ref{thm:harmony=2 semi-harmony} $\big(m(Q),I-Q\big)$ is not harmonious. However, from the proof of Theorem~\ref{thm:matched pair semi harmony is harmony} we see that $\big(m(Q),I-Q\big)$ is semi-harmonious.
\end{example}

We provide another example as follows.

\begin{example}\label{ex:semi-harmonious not implication}
Let $\mathfrak{B}=C(\Omega)$ with $\Omega=[-1,0]\cup [1,2]$ and $\mathfrak{A}=M_2(\mathfrak{B})$, which serves as the Hilbert module over itself. Let $P$ and $Q$ be defined by \eqref{equ:nnew defn of P Q unital Cstar alg} such that
the element $A$ in $\mathfrak{B}$ is specified as
$$A(t)=t,\quad \forall t\in \Omega.$$
Let $\ell$ be the function defined on $\Omega$ by $\ell(t)=\sqrt{t^2-t}$.  It is clear that
\begin{equation*}\ell(A)(t)=\ell(t), \quad \forall\,t\in \Omega, \end{equation*}
where $\ell(A)$ is defined by \eqref{defn of operator ell A}.
Therefore,
\begin{equation}\label{equ:2 l are the same}\ell(A)=\ell(\cdot)\in \mathfrak{B}.\end{equation}
Let $T_1,T_2, T_3,T_4$ and $\mathscr{M}$ be defined by
 \eqref{definition of T1 and T2}, \eqref{equ:new defn of T3 and T4} and \eqref{equ:defn of two orthogonal parts}, respectively. According to Theorem~\ref{thm:invariant submodule M},
$(P|_\mathscr{M},Q|_\mathscr{M})$ is semi-harmonious. In what follows,  we  prove that $(P|_\mathscr{M},Q|_\mathscr{M})$ is not harmonious.

Denote by $b=I_{\mathfrak{B}}-A$. Then
$$b(t)=1-t, \quad \forall t\in\Omega.$$
Let $\mathscr{Y}, \mathscr{Z}\subseteq \mathfrak{B}$ be defined by
\begin{equation*}\mathscr{Y}=\overline{\mathcal{R}\big(b(\cdot)\big)},\quad \mathscr{Z}=\overline{\mathcal{R}\big(\ell(\cdot)\big)}
\end{equation*}
in the sense that for each  $x\in \mathfrak{B}$,
$$\mathcal{R}\big(x(\cdot)\big)=\left\{xy: y\in \mathfrak{B}\right\}\subseteq \mathfrak{B}.$$
Observe that
$$\mathscr{Z}=\overline{\mathcal{R}\big(\ell^2(\cdot)\big)}\subseteq \mathscr{Y},$$ and
$b=bI_{\mathfrak{B}}\in \mathscr{Y}$ with $b(0)=1$, whereas $z(0)=0$ for any $z\in \mathscr{Z}$. It follows that
$\mathscr{Z}\subsetneq \mathscr{Y}$.

With the same convention, by \eqref{equ:2 l are the same} we have
\begin{equation}\label{equ: nullspace-zero-1}\mathcal{N}\big[\ell(A)(\cdot)\big]=\mathcal{N}\big[\ell(\cdot)\big]=\{x\in \mathfrak{B}: \ell(t)x(t)\equiv 0, t\in\Omega\}=\{0\}.
\end{equation}
Following the notations employed in Example~\ref{ex:change to harmony},  we have $\mathscr{M}=\overline{\mathcal{R}(C)}+\overline{\mathcal{R}(D)}$, where
\begin{equation*}C=P(I-Q)P=\left(
          \begin{array}{cc}
            b & 0 \\
            0 &  0\\
          \end{array}
        \right),\quad  D=(I-P)Q(I-P)=\left(
          \begin{array}{cc}
            0 & 0 \\
            0 & b \\
          \end{array}
        \right).
\end{equation*}
So, for any $x\in \mathfrak{A}$ with $x(t)=\big(x_{ij}(t)\big)_{1\le i,j\le 2}$,
\begin{equation}\label{equ:characterization of x in Y}x\in \mathscr{M}\Longleftrightarrow x_{ij}\in \mathscr{Y}\quad\mbox{for $i,j=1,2$}.
\end{equation}
As shown in Example~\ref{ex:change to harmony}, $\overline{\mathcal{R}(T_4|_\mathscr{M})}$ is orthogonally complemented in $\mathscr{M}$ if and only if
\eqref{equ:new orth decomp of M} is satisfied. So, our aim here is to show that
 \begin{equation}\label{equ:strict inclusion-01}\overline{\mathcal{R}(T_4)}+\mathcal{N}(T_3)\cap \mathscr{M}\ne \mathscr{M}.
\end{equation}
For any $x\in \mathfrak{A}$ with $x(t)$ represented as above, by \eqref{nnew expression for T3 T4} and \eqref{equ:2 l are the same}--\eqref{equ:characterization of x in Y}, we conclude that
\begin{align*}&x\in \mathcal{N}(T_3)\cap \mathscr{M}\Longleftrightarrow x_{11},x_{12}\in \mathscr{Y}; x_{21}(t)\equiv 0,x_{22}(t)\equiv 0,\\
&x\in \overline{\mathcal{R}(T_4)}\Longleftrightarrow  x_{11}(t)\equiv 0,x_{12}(t)\equiv 0; x_{21},x_{22}\in \mathscr{Z}.
\end{align*}
Hence, \eqref{equ:strict inclusion-01} can be rephrased as $\mathscr{Z}\ne \mathscr{Y}$, which has already been verified.
\end{example}

\vspace{2ex}

\noindent\textbf{Acknowledgements}\ The authors sincerely thank the referee for his valuable comments and suggestions, which greatly improved the paper's presentation.

\vspace{2ex}

\noindent\textbf{Data Availability}  No data was used for the research described in the paper.

\vspace{2ex}

\noindent\textbf{{\large Declarations}}

\vspace{2ex}

\noindent\textbf{Conflict of interest} No known competing financial interests or personal
relationships that could have appeared to influence the work reported in this paper.

\end{document}